\documentclass{amsart}
\usepackage{amsmath, amsthm}
\usepackage{amsbsy}
\usepackage{amssymb}

\def\be{\begin{equation}}
\def\ee{\end{equation}}

\def\bea{\begin{eqnarray}}
\def\eea{\end{eqnarray}}

\title{Few remarks on the mass spectrum of two-dimensional Toda lattice of
of $E_8$ type}

\author{A. M. Perelomov}

\address{Institute of Theoretical and Experimental Physics,\\  117259 Moscow, Russia}

\begin{document}

\maketitle 

\begin{abstract}  
In this note the simple procedure for obtaining the mass spectrum of two-dimensional 
Toda lattice of $E_8$ type is given.
\end{abstract}

\def\l{\lambda}
\def\R{{\Bbb R}}
\def\bchi{\boldsymbol{\chi}}
\def\ba{\boldsymbol{a}}
\def\pder#1#2{\frac{\partial #1}{\partial #2}}
\def\ppder#1#2#3{\frac{\partial^2 #1}{\partial #2\partial #3}}
\def\ve{\varepsilon}
\def\a{\alpha}
\def\k{\kappa}
\def\l{\lambda}
\def\nn{\nonumber}
\def\ni{\noindent}
\def\bi{\bibitem}
\def\b{\hfil\break}

\section{Introduction. Basic notations} 

Two-dimensional Toda lattice is two-dimensional relativistic field theory
describing $l$ interacting scalar fields. In paper [MOP 1981] it was generalized
for the case of arbitrary simple Lie algebra ${\mathcal G}$ and it has remarkable integrability properties.

This is the relativistic system with Lagrangian

\be
L={1\over 2}\,{\partial _{\mu}} \, {\partial ^{\mu}\, \phi \,} 
- U(\phi),\, \,\, \mu=0,\, 1\,,
\ee
where $\phi = \phi (x_0,\, x_1)$ is $l$-dimensional vector.

The potential $U(\phi )$\, is constructed using some finite set of vectors 
$\{ \alpha _{j} \}$,\\
${j = 0,\, 1,\, ...,\, l}$ in $l$-dimensional Euclidean space 
related to the simple Lie algebra ${\mathcal G}$ of rank $l$ :

\be
U(\phi ) = \sum_{j = 0}^{l} \exp(2\, \alpha _{j},\, \phi ).
\ee

Let us give
some formal definitions, more details may be found in the book [OV 1990].

Let $\mathcal G$ be a compact simple Lie algebra of rank $l$, $R_{+}, (R_{-})$ 
be the set of positive (negative) roots, and $\{\alpha _1,\ldots ,\alpha _l\}$ 
be the set of simple roots. 
Let also $W$ be the Weyl group of root system 
acting in the space $V={{\Bbb R}^l}$, $(\ ,\ )$ be the $W$-invariant bilinear form 
in $V$, $\delta =\sum _{j=1}^{l} n_{j} \alpha _j$ be the highest root, 
$\alpha _{0} = - \delta$, $h = {\sum_{j = 1}^{l} n _{j}} + 1$ 
be the Coxeter number.

In [MOP 1981] the mass spectrum of of scalar fields was found for all 
simple Lie algebras except the most complicated case ${\mathcal G} = E_8$.  
For this algebra only numerical result was given.

\ni In this note we describe two simple methods to obtain  the mass spectrum for the $E_8$ case.
Note that both methods are valid also for the cases of other simple Lie algebras.

\ni
The enumeration of simple roots of the Lie algebra $E_8$ is given on Dynkin diagram.

\begin{center}
\begin{picture}(70,48)(2,-8)
\put(0,0){\circle{8}}
\put(4,0){\line(1,0){30}}
\put(38,0){\circle{8}}
\put(42,0){\line(1,0){30}}
\put(76,0){\circle{8}}
\put(76,4){\line(0,1){30}}
\put(76,38){\circle{8}}
\put(-4,0){\line(-1,0){30}}
\put(-38,0){\circle{8}}
\put(-42,0){\line(-1,0){30}}
\put(-76,0){\circle{8}}
\put(80,0){\line(1,0){30}}
\put(114,0){\circle{8}}
\put(118,0){\line(1,0){30}}
\put(152,0){\circle{8}}
\put(-90,-10){$\alpha_1$}
\put(-52,-10){$\alpha_2$}
\put(-14,-10){$\alpha_3$}
\put(24,-10){$\alpha_4$}
\put(61,-10){$\alpha_5$}
\put(58,36){$\alpha_8$}
\put(99,-10){$\alpha_6$}
\put(137,-10){$\alpha_7$}
\end{picture}
\vspace{3mm}
\footnotesize

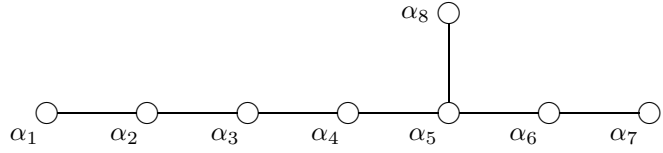
The Dynkin diagram for the Lie algebra $E_8$.
\normalsize
\end{center}

For this enumeration the highest root $\delta $ has the form :
\be \delta =2\alpha _1+3\alpha _2+4\alpha _3+5\alpha _4+6\alpha _5+4\alpha _6+2\alpha _7+
3\alpha _8 . \ee

Note then that in 1989 A.B. Zamolodchikov using conformal theory discover 
that such system appears also at consideration of Ising model in nonzero magnetic field
and he calculated mass spectrum explicitely [Za 1989].

It appears that four mass ratios are equal to the "golden ratio" 
\be
r=\frac{\sqrt{5}+1}{2}=2 \cos({\pi \over{5}}) = 1,6180339887.....
\ee
This remarkable property is related to the fact that for Lie algebra $E_8$ the Coxeter number $h=30$ has factor 5.

In 2010 this theory was confirmed experimentally 
for quasi – one - dimensional Ising ferromagnet
(cobalt niobate) near its critical point [Co 2010]. 

\section{Method 1}
\setcounter{equation}{0}

As it was shown in papers [BCDS 1990], [Fr 1991] masses of paricles are proprtional 
to the components  of special eigenvector of matrix $\rm A$ -- 
Perron - Frobenius vector ( [Pe 1907],\,[Fr 1912] ).  
Here $\rm {A} = 2\,\rm {I} -\rm {C},\, \rm {C} $ is the Cartan matrix  
of Lie algebra $\mathcal G$. For the case $E_8$ we have

\be
{\rm A}\,=\left(
\begin{array}{cccccccc}
0&1&0&0&0&0&0&0\\
1&0&1&0&0&0&0&0\\
0&1&0&1&0&0&0&0\\
0&0&1&0&1&0&0&0\\
0&0&0&1&0&1&0&1\\
0&0&0&0&1&0&1&0\\
0&0&0&0&0&1&0&0\\
0&0&0&0&1&0&0&0
\end{array}
\right) .
\ee

The characteristic equation of this matrix is
\be
x^8 -7 x^6 +14 x^4 -8 x^2 +1\, =\, 0, 
\ee
and his roots are
\be
x_j\, =\, 2\,  \cos(a_j\, \theta),
\ee
where
\be
\theta \, =\, {\pi \over h},
\ee
$h\, =\, 30$\, is Coxeter number.

The numbers 
$\{ a_j\}\, =\, (\,1,\,7,\,11,\,13,\,17,\,19,\,23,\,29)$
are so called exponents of Lie algebra $E_8$.

Note that they have not common divisors with Coxeter number $30$.

Note also that 
$x_5 =- x_4,\,x_6 =- x_3,\,x_7 =- x_2,\,x_8 =- x_1 $ 
and let us give the expressions for $x_j$ in terms of radicals.

\be 
\begin{split}
x_1 = {1\over 2}\,\sqrt{(7+ \sqrt{5}+ \sqrt{30+6 \sqrt{5}}} ,  
x_2 = {1\over 2}\,\sqrt{(7+ \sqrt{5}- \sqrt{30+6 \sqrt{5}}} , \\ 
x_3 = {1\over 2}\,\sqrt{(7- \sqrt{5}+ \sqrt{30-6 \sqrt{5}}} , 
x_4 = {1\over 2}\,\sqrt{(7- \sqrt{5}- \sqrt{30-6 \sqrt{5}}}. 
\end{split}
\ee

The matrix $\rm A$ has nonnegative elements and according to 
Perron - Frobenius theorem [Pe 1907],\,[Fr 1912] it has unique eigenvector
\[
u = (u_1,u_2,u_3,u_4,u_5,u_6,u_7,u_8)
\]
all components of which are positive.
It corresponds to the maximal eigenvalue $\lambda = 2 \cos( \theta )$ 
and we have
\be
u\, {\rm A}\, =\, \lambda\, u\, , 
\ee
\be
\begin{split} 
u_2 = \lambda\, u_1 ,\, u_1+u_3 = \lambda\,\ u_2,\, u_2+u_4 = \lambda\,\ u_3,\, 
u_3+u_5 = \lambda\, u_4,\, \\
u_4+u_6+u_8 = \lambda\, u_5,
u_5+u_7 = \lambda\, u_6,\, u_6 = \lambda\, u_7, 
u_5 = \lambda\, u_8\,.
\end{split}
\ee

Solving these equations, fixing  $u_1 = 2  \sin(\, \theta \,)$, we obtain:

\be
u=\left(
2  \sin(\, \theta \,),
2  \sin(2 \, \theta \,), 
2  \sin(3 \, \theta \,), 
2  \sin(4 \, \theta \,), 
2  \sin(5 \, \theta \,), \\
{{\sin(2 \, \theta \,)} \over  {\sin(3 \, \theta \,)}}, 
{{\sin( \, \theta \,)} \over  {\sin(3 \, \theta \,)}}, 
{{\sin(\, \theta \,)} \over {\sin(2 \, \theta \,)}}\right). 
\ee

\be
u = ( 0,2091;\, 0,4158;\, 0,6180;\, 0,8135;\, 1;\, 0,6728;\, 0,3383;\, 0,5028\, ).
\ee

Note that from this it follows relationes

\be
 {u_7 \over u_1} =r,\, {u_6 \over u_2} = r,\, {u_5 \over u_3} = r,\, {u_4 \over u_8} = r, 
\ee 
where
\be
r=\left({1 \over 2}\right)(1+\sqrt{5})\,\,=2 \cos \left({\pi \over 5}\right)
\ee
is the so called golden section.

This is very nice solution, because these expressions for $u_j$ may be written 
immediately just looking to the Dynkin diagram for $E_8$.

We would like to underline that the situation for arbitrary 
simple Lie algebra  is the same, i. e. the solution may be written 
just to looking to corresponding Dynkin diagram.

Let us give also the expressions for some trigonometric quantities 
in terms of radicals

\be 
\begin{split}
2 \cos({\pi \over{5}}) = {1 \over 2} (1+\sqrt{5}) = r,\,\,
2 \sin({\pi \over{5}}) = \sqrt{{5-\sqrt{5}}\over{2}},\, \quad\quad\quad\quad \quad \quad \quad\quad\\
2 \cos({\pi \over{10}}) = \sqrt{{5+\sqrt{5}}\over{2}},\,
2 \sin({\pi \over{10}}) = \sqrt{{3-\sqrt{5}}\over{2}}, \quad\quad\quad\quad \quad \quad \quad\quad\quad\quad\\
2 \cos({\pi \over{15}}) = {1 \over 2} \sqrt{9 + \sqrt{5} 
+ 2 \sqrt{3} {\sqrt{{5- \sqrt{5}} \over{2}}}},\,
2 \sin({\pi \over{15}}) = {1 \over 2} \sqrt{7 - \sqrt{5} 
- 2 \sqrt{3} {\sqrt{{5- \sqrt{5}} \over{2}}}},\\
2 \cos({\pi \over{30}}) = {1 \over 2} \sqrt{7 + \sqrt{5},\,
+ 2 \sqrt{3} {\sqrt{{5 + \sqrt{5}} \over{2}}}},\,
2 \sin({\pi \over{30}}) = {1 \over 2} \sqrt{9 - \sqrt{5} 
- 2 \sqrt{3} {\sqrt{{5 + \sqrt{5}} \over{2}}}}
\end{split}
\ee

\section{Method 2}
\setcounter{equation}{0}
In the paper [MOP 1981] it was shown that squares of masses are eigenvalues of
matrix
\be
B_{a , b} = \sum _{j=0}^{l} {n_j} {\alpha _j}^a {\alpha _j}^b\,\, n_0\,=\, 1  
\ee
where quantities ${\alpha _j}^a$ are components of vector $\alpha _j$.

For the case of Lie algebra $E_8$ characteristic polynomial $P$ of this matrix is
\be
P = x^8 - 60 x^7 + 1440 x^6 - 18000 x^5 + 127440 x^4
-518400 x^3 + 1166400 x^2-1296000 x + 518400.
\ee
 
In paper [BCDS 1990] was noted factorisation of the characteristic polynomial.

\be
P = P_1\, P_2,\,P_1 = x^4 - 30 x^3 + 240 x^2 - 720 x + 720,\,
P_2 = x^4 - 30 x^3 + 300 x^2 - 1080 x + 720.
\ee

It is easy to check that the roots of polynomial $P_1$ are 
\be
{m_1}^2,\,{m_3}^2,\,{m_4}^2,\,{m_6}^2 
\ee
and the roots of polynomial $P_2$ are
\be
{m_2}^2,\,{m_5}^2,\,{m_7}^2,\,{m_8}^2 .
\ee

Note that 
\be
u_2 u_5 u_7 u_8 = u_1 u_3 u_4 u_6 ,
\ee

and
\be
{m_j}^2 = M {u_j}^2 .
\ee

The quantity $M$ may be found from the equation

\be
M^4 ( u_2 u_5 u_7 u_8 ) ^2 = 720 
\ee

and it has the form

\be
M = 2 \sqrt{3}\,\, {{\sin(6 \theta)} \over {\sin(\theta)}} .
\ee

So the formulae (3.7) and (3.9) give the relation between methods 1 and 2.

Let us give also the explicit expression for 
quantities ${m_j}^2$ in terms of radicals:
\be
\begin{split}
{m_5}^2\, =\, {1 \over 2} \sqrt{15 + 3\,\sqrt{5} 
+ \sqrt{6} {\sqrt{{25 + 11\,\sqrt{5}}}}}\,,\\
{m_7}^2\, =\, {1 \over 2} \sqrt{15 + 3\,\sqrt{5} 
- \sqrt{6} {\sqrt{{25 + 11\,\sqrt{5}}}}}\,,\\
{m_8}^2\, =\, {1 \over 2} \sqrt{15 - 3\,\sqrt{5} 
+ \sqrt{6} {\sqrt{{25 - 11\,\sqrt{5}}}}}\,,\\
{m_2}^2\, =\, {1 \over 2} \sqrt{15 - 3\,\sqrt{5} 
- \sqrt{6} {\sqrt{{25 - 11\,\sqrt{5}}}}}\,,\\
{m_4}^2\, =\, {1 \over 2} \sqrt{15 + 3\,\sqrt{5} 
+ \sqrt{6} {\sqrt{{5 -\,\sqrt{5}}}}}\,,\\
{m_6}^2\, =\, {1 \over 2} \sqrt{15 + 3\,\sqrt{5} 
- \sqrt{6} {\sqrt{{5 -\,\sqrt{5}}}}}\,,\\
{m_3}^2\, =\, {1 \over 2} \sqrt{15 - 3\,\sqrt{5} 
+ \sqrt{6} {\sqrt{{5 +\,\sqrt{5}}}}}\,,\\
{m_1}^2\, =\, {1 \over 2} \sqrt{15 - 3\,\sqrt{5} 
- \sqrt{6} {\sqrt{{5 +\,\sqrt{5}}}}}\,.
\end{split}
\ee

\section{Conclusion}
\setcounter{equation}{0}

The remarkable property of system under consideration is that four mass ratios  in (2.10) are equal to 
the "golden ratio" ("golden section")
\be
r=\frac{\sqrt{5}+1}{2}=2 \cos({\pi \over{5}}) = 1,6180339887.....
\ee
This is only one more fenomenon in which golden ratio appeared.
Note that golden ratio has very long history, see for example book [Co 1961], Ch. 11.
The first book on this topic  "Divina Proportione", illustrated by Leonardo da Vinci, 
was published by Italian mathematician Luca Paccioli in 1509 [Pa 1509].
In conclusion I would like  to give here the quotation of
outstanding astronomer and mathematicien Johannes Kepler [Ke 1596]:
"Geometry has two treasures: one of them is the Pythagorean theorem,
and the other is dividing the segment in average and extreme respect ...
The first can be compared to the measure of gold; the second is more like a gem".

\end{document}